\documentclass{article}
\usepackage{amsmath,amssymb,amsthm}
\usepackage{natbib}
\usepackage{graphicx}

\begin{document}
\DeclareGraphicsExtensions{.eps}

\def\Prodi{\text{\lower0.35cm\hbox{\includegraphics{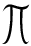}}}}

\title{Causal Inference for Complex Longitudinal Data:\\
The Continuous Time $g$-Computation Formula}

\author{Richard D. Gill
\\
\\
Mathematical Institute, University of Utrecht, Netherlands
\\
EURANDOM, Eindhoven, Netherlands}

\date{November 1, 2001}

\maketitle

\begin{abstract}
I write out and discuss how one might try to prove
the continuous time $g$-computation formula, in 
the simplest possible case: treatments 
(labelled $a$, for actions) and covariates 
($l$: longitudinal data) form together a bivariate counting process.
\end{abstract}

\section{Introduction}

\citet{robins97} outlines a theory of causal inference for complex 
longitudinal data, when treatments can be administered and covariates
observed, continuously in time. This theory is supposed to extend the
earlier work of \citet{robins86,robins87,robins89,robins97}, devoted to 
the case in which covariates and treatments take values in discrete 
spaces, and time advances in discrete time steps. Already in 
\citet{gillrobins01}, we managed to extend the theory to continuously 
distributed covariates and treatments. In this note, we address the 
generalization to continuous time. The major part of this research 
programme has 
already been carried out by \citet{lok01,lok04}. It is an open problem to complete 
that project with a 
continuous time version of the $g$-computation formula and the 
theorems centered around it. The formula tells one how to write down the 
probability distribution of an outcome of interest, in the 
counterfactual situation that a prechosen treatment plan $g$ had been adhered 
to, rather than the factual case that treatment was assigned haphazardly.

\citet{lok01} manages to develop a martingale and counting process 
based theory of \citeauthor{robins97}' (\citeyear{robins97}) 
statistical models, estimators and tests, without having recourse to 
the $g$-computation formula. So is it so central to the theory, after 
all? The answer is that without the formula, the statistical 
methodology lacks motivation. In particular, one needs the formula in 
order to show that the test statistics of \citet{lok01} really do 
test the null hypothesis of no treatment effect, in the sense that 
the counterfactual outcome under all treatment plans $g$ has exactly 
the same probability distribution.

Below we do not succeed in proving the formula, nor establishing the 
wished-for results which should follow from it. What we do do, is present a 
framework in which these questions can hopefully be studied, and in particular, 
write down a conjectural $g$-computation formula and the assumptions 
under which it is likely to be true.

\section{The model}

Suppose that as a patient is followed in time, longitudinal data is 
gathered and treatment decisions or actions are taken; both 
continuously in time. The most simple possible of scenarios, is that 
there is only one kind of action. The only variation in treatment 
is in the times at which the action is taken, the nature of the 
actions at different times is irrelevant or always the same; similarly, 
incoming data  takes the form of a sequence of events at random time points, 
and the only relevant thing is the time of the events, not their nature.
Finally we suppose that actions and longitudinal data events are never
simultaneous. The pair of point processes therefore forms a bivariate 
counting process $(\mathbf N^{a},\mathbf N^{l})$; or if you prefer, a 
single marked point process $\boldsymbol\mu$
with a mark space $\mathcal X=\{a,l\}$, say,
and component point processes $\boldsymbol\mu^{a}$, $\boldsymbol\mu^{l}$; or if you 
prefer, two sequences of random positive time points with no ties between them,
$(0<T^{a}_{1}<T^{a}_{2},\ldots)$, $(0<T^{l}_{1}<T^{l}_{2},\ldots)$.
Ordinary random variables are set in plain lettertype, random processes and 
random measures in bold.
We suppose time varies through a bounded time interval 
$\mathcal T =[0,\tau]$ and that the total number of events of both types is 
finite with probability $1$. Recall that a marked point process is 
a random measure assigning mass $1$ to random ordered pairs  
of a timepoint and accompanying mark, while a counting process counts numbers 
of events, of each kind, up to each timepoint. We suppose there is no event at 
time zero.
The relations between these quantities 
are: $\boldsymbol\mu=\sum_{j}\delta_{(T^{a}_{j},a)}+
\sum_{k}\delta_{(T^{l}_{k},l)}$ where $\delta_{(t,x)}$ is the measure 
with point mass $1$ at the point $(t,x)\in \mathcal T\times\mathcal X$;
$\boldsymbol\mu^{x}(B)=\mu(B\times\{x\})$ for each Borel set in $\mathcal T$
and each mark $x=a,l\in\mathcal X$; $\mathbf N^{x}(t)=\mu^{x}([0,t])$ for 
each $x\in\mathcal X$.

We suppose that we have access to unlimited observational data, and 
therefore essentially know the probabilility distribution, for a randomly 
chosen patient, of the just introduced random quantities. The 
probability law can be recovered from the cumulative intensity process 
or compensator $\boldsymbol\Lambda$ of the 
counting process $\mathbf N$ or, if you prefer, the dual predictable projection or 
compensator $\boldsymbol\nu$ of 
the marked point process $\boldsymbol\mu$. Let $\mu$ (plain lettertype) denote a 
possible realization of the random point process $\boldsymbol\mu$ 
(bold).  Write $\mu_{t}$ for the restriction of the measure 
$\mu$ to $[0,t]\times\mathcal X$. Then for each history of the point 
process up to the time of an event, thus for each $\mu_{t}$ for which there is an 
event at timepoint $t$, we have two 
conditional hazard measures 
$\nu^{x}(\cdot \, |\, t,\mu_{t})$ on $(t,\tau]$, $x=a,l$,
such that the conditional probability that the first event of 
$\boldsymbol\mu$ after $t$ 
is in the time 
interval $\mathrm d s$ and has mark equal to $x$ , given the history up 
to and including time $t$, is $\nu^{x}(\mathrm d s)$ for $s\in(t,\tau]$ and 
$x=a,l$. The two conditional hazard measures have no atoms in common, 
since we assumed there are no simultaneous events.
The dual predictable projection of $\boldsymbol\mu$ is the random
measure $\boldsymbol\nu$ 
defined by 
$\boldsymbol\nu(\mathrm d s,\mathrm d x)=\nu^{x}(\mathrm d s\,|\,t,\mu_{t})$
on the event where $t$ is the time of the last event of $\mu$ strictly 
before time $s$. The cumulative intensity process $\boldsymbol\Lambda$ is 
defined by $\boldsymbol\Lambda^{x}(s)=\boldsymbol\nu((0,s]
\times\{x\})$ for all $s$  and $x$. Thus 
$\boldsymbol\Lambda^{x}(\mathrm ds)=\boldsymbol\nu^{x}(\mathrm d s)=
\nu^{x}(\mathrm d s\,|\,t,\mu_{t})$ where $t$ is as before.

One can generate the whole process by drawing subsequent time points 
and marks using the two conditional hazard measures, given any history 
of events up to the $j$th event at some time point $t$, to generate the 
time and mark of the $j+1$st event.

\section{Treatment plans}

A treatment plan $g$ consists of subplans, one for each $j$ and
$t_{0}^{l}=0<t_{1}^{l}<t_{2}^{l}<\ldots<t_{j}^{l}$, which prescribes 
subsequent action timepoints, from time $t_{j}^{l}$ onwards, 
so long as no further longitudinal data timepoint intervenes. 
We may therefore further split the subplans into sub-subplans, one for 
each $j$ and each $k$, which prescribe the time of the $k$th action 
timepoint after the $j$th longitudinal data timepoint, so long as no 
new longitudinal data timepoint occurs.
The moment there is a new longitudinal data timepoint, the old subplan 
(or subsubplan),
is discarded in favour of the relevant new subplan. Each subplan 
``assumes'' that the overall plan $g$ has been adhered to in previous 
segments of the history, so each subplan ``knows'' all the preceding, 
planned, action timepoints as well as the given preceding longitudinal data 
timepoints. Thus, if we are adhering to a particular 
plan $g$, we can for any sequence of longitudinal data timepoints
$t_{0}^{l}=0<t_{1}^{l}<t_{2}^{l}<\ldots$, thus for any outcome 
$\mu^{l}$, write down the complete 
accompanying sequence of planned action timepoints, and thereby 
reconstruct a complete outcome of a marked point process $\mu^{g}$ 
given the component marked point process outcome $\mu^{l}$. Moreover this can be done in an 
adaptive way: $\mu^{g}_{t}=\mu^{g}|_{(0,t]}$ is a function of $\mu^{l}_{t}
=\mu^{l}|_{(0,t]}$, and of course of the specific treatment plan $g$ under 
consideration.
We can therefore also compute, in an adaptive way, an outcome 
$\Lambda^{g}=(\Lambda^{g,a},\Lambda^{g,l})$  of the 
cumulative intensity process $\boldsymbol\Lambda$, or an outcome 
$\nu^{g}=(\nu^{g,a},\nu^{g,l})$ of the dual predictable 
projection $\boldsymbol\nu$, through its dependence on
$\mu=\mu^{g}$, as a function of any sequence of longitudinal data timepoints
$t_{0}^{l}=0<t_{1}^{l}<t_{2}^{l}<\ldots$, i.e., as a function of 
$\mu^{l}$.

\section{g-Computation Formula}

Suppose a plan $g$ is given. Suppose moreover is given, a random 
variable $Y$, taking values in some Polish space, which we consider 
as the outcome of interest. Alongside the ``factual'' outcome $Y$ we 
suppose there is also defined the ``counterfactual'' outcome $Y^{g}$: 
the outcome which would have pertained, had plan $g$ been adhered to.
Now the conditional law of $Y$ given $\boldsymbol\mu$ can be considered as a 
function of $\mu$, as such we denote it as 
$\mathrm{Law}(Y|\boldsymbol\mu=\mu)$. Therefore, 
for a given sequence of longitudinal data timepoints
$t_{0}^{l}=0<t_{1}^{l}<t_{2}^{l}<\ldots$, which determines a possible 
outcome of $\mu^{l}$, we can evaluate the law of 
$Y$ given $\boldsymbol\mu$ at 
$\boldsymbol\mu=\mu^{g}=\mu^{g}(\mu^{l})=(\mu^{l},\mu^{a}(\mu^{l},g))$.
The $g$-computation formula, which we want to prove under versions of 
the usual three assumptions of consistency, no-unmeasured confounding, 
and evaluability, is the following:
\begin{equation}
    \begin{aligned}
 \    \mathrm{Law}&(Y^{g})~=~
\sum_{n}\idotsint\limits_{t^{l}_{1}<\ldots<t^{l}_{n}\le\tau} \\
    &\prod_{i=1}^{n}\Prodi_{s\in 
    (t^{l}_{i-1},t^{l}_{i})}\Bigl(1-\Lambda^{g,l}(\mathrm d s)\Bigr)
          \Lambda^{g,l}(\mathrm d t^{l}_{i})
             \Prodi_{s\in 
    (t^{l}_{n},\tau]}\Bigl(1-\Lambda^{g,l}(\mathrm d s)\Bigr)
    \mathrm{Law}(Y|\boldsymbol\mu=\mu^{g}).
    \end{aligned}
    \notag
\end{equation}	 
The first thing to note about this formula is that it is a functional 
of the cumulative intensity function $\Lambda^{g,l}$ and of the 
conditional law of $Y^{g}$ given $\boldsymbol\mu$, both considered as 
functionals of $\mu^{g}$, which again is a functional of the chosen 
treatment plan $g$ and the summation and integration variables in the 
formula: the 
total number $n$ of longitudinal data timepoints in the time interval $\mathcal T$
and their values $0=t^{l}_{0}<t^{l}_{1}<\ldots<t^{l}_{n}\le\tau$. 
These variables precisely determine an outcome of $\boldsymbol\mu^{l}$.
The cumulative intensity function  $\Lambda^{g,l}$ is computed from 
the conditional probability laws of the `next longitudinal data 
timepoint' restricted to the event, that it precedes the next action timepoint,
given the history of the process $\boldsymbol\mu$ up to the times of the 
zero'th, first, second \ldots\ events. Thus it depends on which 
version is chosen of each of these conditional probability laws. 

Recall from \citet{gillrobins01} that there are two issues in establishing this formula. The 
first is the question whether, when one chooses appropriate versions 
of the conditional distributions involved, it gives the right answer. 
The second question is whether, when conditional distributions are 
chosen, if possible, in some canonical fashion, the result 
is uniquely defined as a functional of the joint law of the data 
$\boldsymbol\mu,Y$. We may have to face up to one third, more technical issue: 
the formula supposes that in the counterfactual world where treatment 
plan $g$ is followed, there is no explosion in the sequence of 
timepoints of events; in other words, if we replace the conditional 
law of $Y$ in the integrand with the constant function $1$, the 
result of the g-computation formula should be the total probability $1$.
Let us call this condition, the no-explosion condition for plan $g$.

Now we discuss what the three usual conditions should look like, in 
this context, and make some remarks on how one might attempt to prove 
the formula.

The consistency condition, in a sufficient and weaker `in law' form, 
should naturally be: 
$\mathrm{Law}(Y|\boldsymbol\mu=\mu)=\mathrm{Law}(Y^{g}|\boldsymbol\mu=\mu)$
for outcomes $\mu$ consistent with plan $g$: thus, outcomes $\mu$ such that 
$\mu^{a}=\mu^{a}(\mu^{l},g)$. The `no unmeasured confounders' 
assumption should be that the intensity process of the action events, 
when the history of the process $\boldsymbol\mu$ is augmented by taking $Y^{g}$ 
to be a random variable realized at time $t=0$, should be the same as 
the intensity process of the action events when only the history of 
$\boldsymbol\mu$ is taken into account, 
for outcomes $\mu$ consistent with plan $g$. 
In terms of conditional distributions, it is the 
assumption that conditional on the times and types of events up to 
any number of the events, $Y^{g}$ is independent of the time to the next 
action event, restricted to the event that it precedes the next 
longitudinal data event; and we only need to check this condition for 
outcomes $\mu$ consistent with plan $g$.  Just the consistency and the no unmeasured 
confounders assumptions should be sufficient to establish the 
correctness of the $g$-computation formula, when the same conditional 
distributions are employed in the formula, as are involved in the 
assumptions. Since typically the probability that $\boldsymbol\mu$ is consistent 
with $g$ is zero, this result has no empirical content. Still, given 
versions of all involved conditional distributions, the result is 
not obviously true, so does have mathematical content. 
The first step in the proof is naturally to 
replace $Y$ with $Y^{g}$ on the right hand side of the formula, 
using the consistency assumption. How to proceed from here, is not so 
clear. A strategy which might work, is to consider the right hand 
side of the $g$-computation formula, with $Y$ replaced by $Y^{g}$ 
and $\tau$ replaced by a variable timepoint $\sigma\in\mathcal T$ as 
a function of $\sigma$, say $b(\sigma)$, and show that it satisfies some integral 
equation. We are given the value of the function $b$ at $\sigma=\tau$. 
If one can show the integral equation is uniquely solved by a constant function 
$b^{*}$ satisfing $b^{*}(0)=\mathrm{Law}(Y^{g})$, we are done. 
The non-explosion condition will presumably be 
needed in this analysis. The important step is guess a non-trivial
probabilistic interpretation of $b(\sigma)$, and 
take the guess to define a function $b^{*}(\sigma)$. Next,
use the probabilistic interpretation to write informally a 
relation between $b^{*}(\sigma+\mathrm d s)$ and 
$b^{*}(\sigma)$, as an expectation of the possible 
outcomes in the time interval $\mathrm d s$. Use probability theory to 
convert this to a rigorous relation in integral form. 

Informally, the proof should parallel that in the discrete time case 
and correspond to the remark that the law of $Y^{g}$ given 
$\mu_{\sigma+\mathrm d \sigma}$ does not depend on 
$\mu^{a}(\mathrm d \sigma)$. Therefore, in order to 
recover the law of $Y^{g}$ given $\mu_{\sigma}$ by averaging over the 
conditional law of the events of $\mu$ in the time interval $\mathrm d \sigma$ 
given the events in the past, we need only average over the 
conditional law of the longitudinal data events. But whether or not 
there is a longitudinal data timepoint in this small time interval is
a Bernoulli $(\Lambda^{l}(\mathrm d \sigma))$ variable. Thus 
$\mathrm{Law}(Y^{g}|\mu_{\sigma})$ is a Bernoulli $(\Lambda^{l}(\mathrm d \sigma))$
mixture of the two distributions
$\mathrm{Law}(Y^{g}|\mu_{\sigma+\mathrm d \sigma})$ with 
$\mu^{l}(\mathrm d \sigma)=0,1$.

Another possible ingredient is yielded by the remark that the law of 
$Y^{g}$ given $\boldsymbol\mu_{t}$ is a martingale in $t$ with 
respect to the history of $\boldsymbol\mu$, and hence can be written 
as a stochastic integral with respect to 
$\boldsymbol\mu-\boldsymbol\nu$. The representation involves the 
intensities of $\boldsymbol\mu$ with respect to its own history, and 
with respect to the augmented history when  $Y^{g}$ is realized at 
time $0$.

In order to obtain a result with empirical content, we have to show 
how the formula can be uniquely evaluated, under further assumptions, 
from the joint law of $Y$ and $\boldsymbol\mu$. A natural assumption which 
guarantees a canonical choice of conditional laws is continuity: 
we should assume that versions of all the conditional laws involved in 
the $g$ computation formula, can be chosen so as to be continuous 
on the support of the conditioning variables. The 
conditioning variables are partial histories of $\boldsymbol\mu$ up to the 
so-manyth event, and the total history of $\boldsymbol\mu$ on $\mathcal T$. 
Continuity of probability laws is in the sense of weak convergence, 
and the partial and total histories of $\boldsymbol\mu$ are given their natural 
topologies. The conditional laws now have canonical versions on the 
supports of the conditioning variables, and we should make the 
evaluability condition on the plan $g$ that 
for partial histories in the support of the corresponding partial 
history of $\boldsymbol\mu$, the next planned action time (restricted to the 
event where it precedes the next longitudinal data timepoint)
lies in the support of the conditional distribution of that time given 
the partial history so far.

\bibliographystyle{Chicago}
\raggedright

\bibliography{cicld-cc}

\end{document}